\documentclass[12pt]{amsart}
\usepackage{amssymb}

\begin{document}
%\vsize=7.85in
%\hsize=14cm
%\oddsidemargin=-1.2cm\oddsidemargin=-1cm
%\evensidemargin=-1.2cm\evensidemargin=-1cm
%\hoffset 5.5truemm
%\setlength{\textwidth}{25cm}
\numberwithin{equation}{section}

\def\Label#1{\label{#1}}

\def\1#1{\overline{#1}}
\def\2#1{\widetilde{#1}}
\def\3#1{\mathcal{#1}}
\def\4#1{\widehat{#1}}
\def\5#1{\mathbb{#1}}

\def\th{\theta}
\def\ra{\rightarrow}
\def\p{\partial}
\def\a{\alpha}
\def\b{\beta}
\def\t{\theta}
\def\g{\gamma}
\def\wt{\widetilde}
\def\O{\Omega}
\def\l{{\ell}}
\def\II{\sqrt{-1}}
\def\e {\epsilon}
\def\PP{{\5 P}}
\def\BB{{\5 B}}
\def\HH{{\5 H}}
\def\db{\overline{\partial}
\def \sm{\setminus}
\def\CC{{\5 C}}
\def\RR{{\bf R}}
\def\ZZ{{\bf Z}}
\def\KK{{\5 K}}
\def\-{\overline}
\def\ale{\mathrel{\mathop{<}\limits_{\sim}}}
\def\o{\omega}
\def\e{\epsilon}
\def\D{\Delta}
\def\cal {\Cal}
\def\ss{\subset}
\def\L{\Lambda}
\def\h{\hbox}
\def\d{\delta}
\def\lll{\tau, \lambda}
\def\dv\{\hbox{d(Vol)}}
\def\du\{\hbox{d(Eucl)}
\def\wh{\widehat}
\def\H{\hat}

\def\d{\delta}
\def\a{\alpha}
\def\M*{\wt{M^*}}
\def\Ml{\wt {M_{\l}}}
\def\C{\5 {C}}
\def\-{\overline}
\def\ale{\mathrel{\mathop{<}\limits_{\sim}}}
\def\age{\mathrel{\mathop{>}\limits_{\sim}}}
\def\ld{\lambda}
\def\O{\Omega}
\def\o{\omega}
\def\ow{o_{wt}}
\def\ss{\subset}
\def\sm{\setminus}
\def\j{j,j^*}
\def\L{\Lambda}
\def\pd{\phi(\Delta)}
\def\h{\hbox}
\def\HX{\hat{X}}
\def\df{{d\Phi^*\OO d\xi}}
\def\q{Q(X,\xi)}
\def\wt{\widetilde}
\def\ra{\rightarrow}
\def\ih{I_{2n+1}+i{\partial H(X)\over\partial X }}
\def\l{\ell}
\def\Sqt{\sqrt}
\def\H{\hat}
\def\a{\alpha}
\def\d{\delta}
\def\V{\nu}
\def\Ol{\overline}
\def\OO{\over}
%\hoffset=-0.1in
%\voffset=1.0in
%\hsize=7.0truein
%\vsize=8.6truein
\def\s{\sigma}
\def\D{\Delta}
\def\d{\delta}
\def\b{\beta}
\def\a{\alpha}
\def\Sqt{\sqrt}
\def\bC{{\5C}}
\def\bR{\5 R}
\def\<{\langle}
\def\>{\rangle}
\def\id{{\sf id\,}}
\def\Aut{{\sf Aut}}
\def\im{\text{{\sf Im\,}}}
\def\Im{\text{{\sf Im\,}}}
\def\re{\text{{\sf Re\,}}}
\def\Re{\text{{\sf Re\,}}}
\def\dbl{[\![}
\def\dbr{]\!]}
\def\R{\text{{\sf R}}}
\def\s{\text{{\sf S}}}
\def\:{\colon}

\title[The equivalence problem and rigidity]
{The equivalence problem and rigidity for hypersurfaces embedded into
hyperquadrics}
\author[P. Ebenfelt, X. Huang, D. Zaitsev]
{Peter Ebenfelt, Xiaojun Huang, and Dmitri Zaitsev}
\footnotetext{{\rm The first author is supported in part by NSF-0100110 and
a Royal Swedish Academy of Sciences Research Fellowship.
The second author is supported in part by NSF-0200689 and a grant from
the Rutgers Research Council. The third author is supported in part
by a grant from the Italian Consiglio Nazionale delle Ricerche and
a research grant of the University of Padua}}
\address{P. Ebenfelt: Department of Mathematics, University
of California at San Diego, La Jolla, CA 92093-0112, USA}
\email{pebenfel@math.ucsd.edu }
\address{X. Huang: Department of Mathematics, Rutgers University at New
Brunswick, NJ 08903, USA} \email{huangx@math.rutgers.edu}
\address{D. Zaitsev: Dipartimento di Matematica, Universit\`a di Padova, via Belzoni
7, 35131 Padova, Italy} \email{zaitsev@math.unipd.it}

%\address{D. Zaitsev: Mathematisches Institut,
%Eberhard-Karls-Universit\"at T\"ubingen, D-72076 T\"ubingen, GERMANY}
%\email{dmitri.zaitsev@uni-tuebingen.de}

\newtheorem{Thm}{Theorem}[section]
\newtheorem{Cor}[Thm]{Corollary}
\newtheorem{Pro}[Thm]{Proposition}
\newtheorem{Lem}[Thm]{Lemma}
\theoremstyle{definition}
\newtheorem{Def}[Thm]{Definition}
\newtheorem{Rem}[Thm]{Remark}
\newtheorem{Exa}[Thm]{Example}

\def\bl{\begin{Lem}}
\def\el{\end{Lem}}
\def\bp{\begin{Pro}}
\def\ep{\end{Pro}}
\def\bt{\begin{Thm}}
\def\et{\end{Thm}}
\def\bc{\begin{Cor}}
\def\ec{\end{Cor}}
\def\bd{\begin{Def}}
\def\ed{\end{Def}}
\def\br{\begin{Rem}}
\def\er{\end{Rem}}
\def\be{\begin{Exa}}
\def\ee{\end{Exa}}
\def\bpf{\begin{proof}}
\def\epf{\end{proof}}
\def\ben{\begin{enumerate}}
\def\een{\end{enumerate}}

%\abstract We consider the class of Levi nondegenerate
%hypersurfaces $M$ in $\bC^{n+1}$ that admit a local (CR
%transversal) embedding, near a point $p\in M$, into a standard
%nondegenerate hyperquadric in $\Bbb C^{N+1}$ with codimension
%$k:=N-n$ small compared to the CR dimension $n$ of $M$. We show
%that, for hypersurfaces in this class, there is a normal form
%(which is closely related to the embedding) such that any local
%equivalence between two hypersurfaces in normal form must be an
%automorphism of the associated tangent hyperquadric. We also show
%that if the signature of $M$ and that of the standard hyperquadric
%in $\bC^{N+1}$ are the same, then the embedding is rigid in the
%sense that any other embedding must be the original embedding
%composed with an automorphism of the quadric.
%\endabstract

\maketitle

\section{Introduction}\Label{s0}

Our main objective in this paper is to study the class of real
hypersurfaces $M\subset \bC^{n+1}$ which admit holomorphic
(or formal) embeddings
%, centered at a point $p\in M$,
into the unit sphere (or, more generally, Levi-nondegenerate
hyperquadrics) in $\bC^{N+1}$ where the codimension $k:=N-n$ is
small compared to $n$. Such hypersurfaces play an important role
e.g.\ in deformation theory of singularities where they arise as
links of singularities (see e.g.\ \cite{BM}). Another source is
complex representations of compact groups, where the orbits are
always embeddable into spheres due to the existence of invariant
scalar products.
% (or, more precisely, $k< n/2$).
One of our main results is a complete normal form for
hypersurfaces in this class with a rather explicit solution to the
equivalence problem in the following form (Theorem~\ref{t0.3}):
{\em Two hypersurfaces in normal form are locally
biholomorphically equivalent if and only if they coincide up to an
automorphism of the associated hyperquadric}. Our normal form here
is different from the classical one by Chern--Moser~\cite{CM}
(which, on the other hand, is valid for the whole class of Levi
nondegenerate hypersurfaces), where, in order to verify
equivalence of two hypersurfaces, one needs to apply a general
automorphism of the associated hyperquadric to one of the
hypersurfaces, possibly loosing its normal form, and then perform
an algebraically complicated procedure of putting the transformed
hypersurface back in normal form. Another advantage of our normal
form, comparing with the classical one, is that it can be directly
produced from an embedding into a hyperquadric and hence does not
need any normalization procedure.

Our second main result is a rigidity property for embeddings into
hyperquadrics (Theorem~\ref{t0.4}), where we show, under a
restriction on the codimension, that any two embeddings into a
hyperquadric of a given hypersurface coincide up to an
automorphism of the hyperquadric.

Before stating our main results more precisely, we introduce some
notation. Let $\mathbb H^{2n+1}_\ell$ denote the standard
Levi-nondegenerate hyperquadric with signature $\ell$ ($0\leq
\ell\leq n$):
\begin{equation}\Label{0.1}
{\5 H}_{\ell}^{2n+1}:= \Big\{ (z,w)\in {\5 C}^n\times \5 C : \Im
w=-\sum_{j=1}^{\l}|z_j|^2+\sum_{j=\l+1}^n|z_j|^2 \Big\}.
\end{equation}
Since the hyperquadric ${\5 H}_{\ell}^{2n+1}$ is clearly
(linearly) equivalent to ${\5 H}_{N-\ell}^{2n+1}$, we may restrict
our attention to $\ell \le \frac {n}{2}$. When $\l=0$, ${\5
H}_{\l}^{2n+1}$ is the Heisenberg hypersurface, also denoted by ${\5
H}^{2n+1}$, which is locally biholomorphically equivalent to the
unit sphere in $\bC^{n+1}$.
 For brevity, we shall
use the notation $\<\cdot,\cdot\>_\l$ for the standard
complex-bilinear scalar product form of signature $\l$ in $\bC^n$:
% i.e.\ for a fixed $\l\leq n/2$
%and $$, we define
\begin{equation}\Label{0.2}
\< a,b\>_\l:=-\sum_{j=1}^\l a_jb_j+\sum_{j=\l+1}^n a_jb_j, \quad
a,b\in \bC^n.
\end{equation}
The dimension $n$ will be clear from the context and we shall not
further burden the notation by indicating also the dependence of
$\<\cdot,\cdot\>_\l$ on $n$. Recall that a smooth ($C^\infty$)
real hypersurface $M$ in ${\5 C}^{n+1}$ is Levi-nondegenerate of
signature $\ell$ (with $\ell\leq n/2$) at $p\in M$ if it can be
locally approximated, at $p$, by a biholomorphic image of ${\5
H}_{\ell}^{2n+1}$ to third order, i.e.\ if there are local
coordinates $(z,w)\in \bC^{n}\times\bC$ vanishing at $p$ such that
$M$ is defined, near $p=(0,0)$, by
\begin{equation}\Label{1.3}
\Im w=\<z,\bar z\>_\l+A(z,\bar z, \Re w),
\end{equation}
where $A(z,\bar z,u)$ is a smooth function which vanishes to third
order at $0$. (In fact, one can assume, after possibly another
local change of coordinates, that the
function $A$ vanishes at least to fourth order at $0$, cf.~\cite{CM}.)
In this paper, we shall consider {\em formal hypersurfaces}
defined by formal power series equations of the form \eqref{1.3}
and --- more generally --- of the (not necessary graph) form
\begin{equation}\Label{0.3}
\Im w=\<z,\bar z\>_\l+A(z,\bar z,w,\bar w),
\end{equation}
where $A(z,\bar z,w,\bar w)$ is a real-valued
formal power series vanishing at least to fourth order.
Our motivation for
considering equations of this form lies in the fact that
they arise naturally in the study of Levi-nondegenerate
hypersurfaces $M\subset \bC^{n+1}$ admitting
embeddings into hyperquadrics;
see Proposition~\ref{p0.2} and section~\ref{s4} for a detailed
discussion. In this paper, we shall only concern ourselves with
the formal study of real hypersurfaces and, hence, we shall
identify smooth functions $A(z,w,\bar z,\bar w)$ with their
formal Taylor series in $(z,w,\bar z,\bar w)$ at $0$.

Given a formal power series $A(Z,\bar Z)$, $Z=(z,w)\in\bC^{n}\times\C$, we
can associate to it
two linear subspaces $V_A\subset\bC\dbl Z\dbr$, $U_A\subset\bC\dbl
\bar Z\dbr$, where $\bC\dbl X\dbr$ denotes the ring of formal
power series in $X$ with complex coefficients, as follows:
$$
V_A:=\text{{\rm span}}_\bC\Big \{\frac{\partial^\alpha
A}{\partial \bar Z^\alpha}(Z,0)\Big\},
\quad
U_A:=\text{{\rm span}}_\bC\Big \{\frac{\partial^\alpha
A}{\partial
 Z^\alpha}(0,\bar Z)\Big\},
$$
where ``$\text{{\rm span}}_\bC$'' stands for the linear span over
$\bC$ and $\alpha$ runs over the set $\5 N^{n+1}$ of all
multi-indices of nonnegative integers. We shall only consider
formal series $A(Z,\bar Z)$ which are {\em real-valued} (i.e.\
which formally satisfy $A(Z,\bar Z)=\-{A(Z,\bar Z})$). For such
power series one has ${V_A}=\-{U_A}$. We define the {\it rank of
$A$}, denoted $\R(A)$, to be the dimension of $V_A$, possibly
infinite. We shall see (Proposition~\ref{p0.2}) that hypersurfaces
admitting (CR transversal, see below)  embeddings into
hyperquadrics can be represented by equations \eqref{0.3} with $A$
being of finite rank. If $\R(A)=:r<\infty$, then a linear algebra
argument shows (cf.\ \cite{We78,D82,We99}) that there are a
nonnegative number $\s(A)=:s\leq \R(A)$, called here the {\it
signature of $A$}, and formal power series $\phi_j\in \bC\dbl
Z\dbr$, $j=1,\ldots, r$, linearly independent (over $\bC$), such
that
\begin{equation}\Label{0.4}
A(Z,\bar Z)=-\sum_{j=1}^{s} |\phi_j(Z)|^2+\sum_{j=s+1}^{r}
|\phi_j(Z)|^2
\end{equation}
(where as usual, in the special cases $s=0$ and $s=r$, the
corresponding void sums in \eqref{0.4} are understood to be zero).
See also \cite{D01} for conditions on $A$ yielding $s=0$.
Moreover, it holds that the span of the $\phi_j$ is equal to $V_A$
(indeed, in the linear algebra argument alluded to above, the
$\phi_j$ are chosen as a "diagonal" basis for $V_A$) and the
collection of vectors $\big(\frac{\partial^\alpha \phi}{\partial
Z^\alpha}(0)\big)_\a$, where $\phi=(\phi_1,\ldots, \phi_r)$ and
$\a$ is as above, spans $\bC^r$. We should point out here,
however, that for an arbitrary collection $(\phi_j)_{j=1}^r$ in
$\bC\dbl Z\dbr$, a number $s\leq r$ and the corresponding series
$A(Z,\bar Z)$ defined by \eqref{0.4}, it only holds, in general,
that $V_A$ is contained in the span of the $\phi_j$. In the case
where $V_A$ is strictly contained in the latter span, however,
there is another representation of $A$ in the form \eqref{0.4}
with another collection of $r'<r$ formal power series.
% $(\tilde \phi_j)_{j=1}^{r'}$ with $r'<r$.

\bd\Label{d0.1} For each nonnegative integer $k$, define the class
$\3 H_k$ to consist of those real-valued formal power series
$A(z,\bar z,w,\bar w)$ in $(z,w)\in \bC^n\times\bC$ satisfying
either of the following equivalent conditions:
\begin{enumerate}
\item[(i)] $\R(A)\leq k$ and no formal power series in $V_A$ has a constant or linear term;
\item[(ii)] $A(z,\bar z,w,\bar w)=\sum_{j=1}^k \phi_j(z,w)\1{\psi_j(z,w)}$
for some formal power series $\phi_j,\psi_j\in \C \dbl z,w \dbr$
having no constant or linear terms.
\end{enumerate}
\ed

Before stating our main results, we give an auxiliary result
relating $\3 H_k$ to the class of hypersurfaces admitting formal
(CR) embeddings into the hyperquadric $\5H^{2N+1}_{\l'}$. Observe
first that, for $\l'\ge n+1$ (and $\l'\leq N/2$, according to our
convention), $\C^{n+1}$ can be linearly embedded into
$\5H^{2N+1}_{\l'}$ via the map
$$Z\mapsto (Z,0,Z,0)\in \C^{n+1}\times\C^{\l'-n-1}\times\C^{n+1}\times\C^{N-n-\l'}.$$
Hence also any hypersurface $M\subset\C^{n+1}$ can be
``trivially'' embedded into $\5H^{2N+1}_{\l'}$. In order to avoid
this kind of embedding, we restrict our attention here to formal
CR embeddings $H\colon (M,0)\to (M',0)$ (where $M'\subset
\bC^{N+1}$ is another hypersurface) such that
$$dH(T_0M)\not\subset T^c_{0}M':=T_0M'\cap iT_0M'.$$
We shall call such embeddings {\em CR transversal}. In particular,
in the case $M'=\5 H^{2N+1}_{\l'}$ considered in present paper, CR
transversality is equivalent to $\partial H_{N+1}(0)\neq 0$. It is
easy to see that any CR embedding of a hypersurface
$M\subset\C^{n+1}$ into $\5H^{2N+1}$ is automatically CR
transversal. More generally, a CR embedding of $(M,0)$ into
$(M',0)$, where $M'$ is Levi-nondegenerate with signature $\l'<n$
at $0$, is CR transversal. Indeed, if it were not, then, as is
well known, the Levi form of $M'$ at $0$ must vanish on the complex
tangent space of $M$ at $0$ and this cannot happen unless $\l'\geq
n$. We have:

\bp\Label{p0.2}
Let $M$ be a formal Levi-nondegenerate hypersurface
in $\bC^{n+1}$  of signature $\l\leq n/2$ at a point $p$.
Then, $M$ admits a formal CR transversal embedding
%$H=(H_1,\ldots, H_{N+1})$,
into $\5 H^{2N+1}_{\l'}$, $\ell'\leq N/2$, if and only if $\l'\geq
\l$ and there are formal coordinates $(z,w)\in \bC^n\times \bC$,
vanishing at $p$, such that $M$ is defined by an equation of the
form \eqref{0.3} with $A\in \3 H_{N-n}$.

More precisely, if $A\in \3 H_{k}$ for some $k$, then the
hypersurface $M$, given by \eqref{0.3}, admits a formal CR
transversal embedding into $\5 H^{2N+1}_{\l'}$ with $N=n+k$ and
$\l'=\min(\l+\s(A),N-\l-\s(A))$. Conversely, if $M$ admits a
formal CR transversal embedding into $\5 H^{2N+1}_{\l'}$, then
there are formal coordinates $(z,w)$ as above and a formal power
series $A(z,\bar z,w,\bar w)$ such that $M$ is defined by
\eqref{0.3} with $\R(A)\leq {N-n}$ and $\s(A)\leq \max(\l'-\l,
N'-l'-l)$. Moreover, if in addition $\l'<n-\l$, then $\s(A)\leq
\l'-\l$ holds. \ep

We shall use the notation $\Aut(\5 H^{2N+1}_{\l},0)$ for the
stability group of $\5 H^{2N+1}_{\l}$ at $0$, i.e.\ for the group of
all local biholomorphisms $(\bC^{N+1},0)\to (\bC^{N+1},0)$ preserving
$\5 H^{2N+1}_\l$. Recall
%(c.f.\ \cite{CM} [CM] or [BER]) that
%$\Aut(\5 H^{2N+1}_\l)$ is explicitly known and
that every $T\in \Aut(\5 H^{2N+1}_\l,0)$ is a linear fractional
transformation of $\bC^{N+1}$ (see e.g.\ \eqref{autos0} below or
\cite{BERbull} for an explicit formula). Our first main result
states that a defining equation of the form \eqref{0.3} with $A\in
\3 H_k$, $k<n/2$, is unique modulo automorphisms of the associated
quadric. More precisely, we have the following.

\bt\Label{t0.3}
Let $M_j$, $j=1,2$, be formal
Levi-nondegenerate hypersurfaces of signature $\l$ in $\bC^{n+1}$
given by
$$
\im w=\<z,\bar z\>_\l+A_j(z,\bar z,w,\bar w)
$$
respectively, where $(z,w)\in \bC^n\times\bC$.
Assume that $A_j\in \3 H_{k_j}$ with $k_1+k_2<n$.
Then, $(M_1,0)$ and $(M_2,0)$ are formally equivalent if and only if there exists
$T\in \Aut(\5 H^{2n+1}_\l,0)$ sending $M_1$ to $M_2$. More
precisely, any formal equivalence sending $(M_1,0)$ to $(M_2,0)$
is the Taylor series of an automorphism
\begin{equation}\Label{autos0}
T(z,w):= \frac
{(\lambda(z+aw)U,\sigma \lambda^{2}w)}
{1-2i\<z,\bar {a}\>_\l -(r+i\<a,\bar{a}\>_{\l})w}
\in \Aut(\5 H^{2n+1}_\l,0)
\end{equation}
for some $\lambda>0$, $r\in \bR$, $a\in \bC^n$, $\sigma=\pm1$
and an invertible $n\times n$ matrix $U$ with
$\<zU,\bar z\bar U\>_\l = \sigma \<z,\bar z\>_\l$
such that
\begin{equation}\Label{eq-rel}
A_1\equiv \sigma\lambda^{-2} |q|^2 \, A_2\circ (T,\bar T),
\end{equation}
where $q=q(z,w)$ is the denominator in \eqref{autos0}.
\et

\br
 It follows immediately from \eqref{eq-rel}
that $\R(A_1)=\R(A_2)$, and either $\s(A_1)=\s(A_2)$
or $\s(A_1)= \R(A_2)-\s(A_2)$.
Thus, the rank $\R(A)$ and the two-point-set $\{\s(A),\R(A)-\s(A)\}$
%$\min(\s(F_1),\R(F_1)-\s(F_1))$
are invariants of $M$ (and not only of the equation \eqref{0.3})
provided $A\in \3 H_k$ for some $k<n/2$.\er

By combining Theorem~\ref{t0.3} with Proposition~\ref{p0.2}, one
can produce examples of formal real hypersurfaces $M$
%of signature $\l$
which
%, for a fixed $\l'$, with $\l\leq \l'\leq n-\l$ and
%$2\le k\le n-2$, and any $N<2n-k$,
do not admit formal CR transversal embeddings into $\5 H^{2N+1}_{\l'}$.
%Indeed, a sufficient
%condition is that $M$ is given in local coordinates by an equation
%of the form \eqref{0.3}, where $A \in \3 H_k$ and
%$\min(\s(A),\R(A)-\s(A))> \l'-\l$. The verification of this (using
%directly Proposition \ref{p0.2} and Theorem \ref{t0.3}) is left to
%the reader.
Let us give an explicit example.
%with $k=2$ and $\l'=\l<n/2$.

\begin{Exa} Let $M\subset\bC^{n+1}$ be defined by
$$
\im w=\<z,\bar z\>_\l+A(z,\bar z,w,\bar w),
$$
where $\l<n/2$ and
$$
A(z,\bar z,w,\bar w)=\re(w^s\-{h(z)}),$$ for some $s\geq 2$ and $h(z)$
a (nontrivial) homogeneous polynomial of degree $\geq 2$. It is
not difficult to see that $\R(A)=2$ and, since $A$ takes both
positive and negative values, $\s(A)=1$. Hence,
%by the remarks above,
in view of Proposition \ref{p0.2} and Theorem \ref{t0.3},
$(M,0)$ cannot be formally embedded
%, at $0$,
into $\Bbb H^{2N+1}_\l$ for any $N$ with $N<2n-2$. On the other
hand, $M$ can be embedded into ${\Bbb H}_{\l+1}^{2(n+2)+1}$, again
by Proposition \ref{p0.2}.
\end{Exa}

Our second main result of this paper is a rigidity result for
formal embeddings into hyperquadrics of the same signature. Since
the signature $\l$ of $M$ is $\leq n/2$ (and then, in particular,
$<n$) it follows that every formal embedding in this case is
automatically CR transversal.

\bt\Label{t0.4} Let $M$ be a formal Levi-nondegenerate
hypersurface in $\bC^{n+1}$ of signature $\l\leq n/2$ at a point
$p$. If $\l=n/2$, then we shall also assume that $(M,p)$ is not
formally equivalent to $(\mathbb H^{2n+1}_{n/2},0)$. Suppose that,
for some integers $k_1\le k_2$, there are two formal embeddings
$$H_1\: (M,p)\to \left(\5 H^{2(n+k_1)+1}_\l,0\right), \quad
H_2\: (M,p)\to \left(\5H^{2(n+k_2)+1}_\l,0\right).$$ If
$k_1+k_2<n$, then there exists an automorphism $T\in \Aut(\5
H^{2(n+k_2)+1}_\l,0)$ such that $$H_2= T\circ L\circ H_1,$$ where
$L\colon \bC^{n+k_1}\times\bC\to \bC^{n+k_2}\times\bC$ denotes the
linear embedding $(z,w)\mapsto (z,0,w)$ which sends $\5
H^{2(n+k_1)+1}_\l$ into $\5 H^{2(n+k_2)+1}_\l$.
%$$(z,w)\in\C^N\times\C \mapsto (z,0,w)\in\C^N\times\C^{N'-N}\times\C. $$
\et

If $(M,p)$ is formally equivalent to $(\mathbb
H^{2n+1}_{n/2},0)$, then the conclusion of Theorem~\ref{t0.4}
fails. Indeed, if $n=2\l$, then the two
embeddings $L, L_-\colon(\mathbb H^{2n+1}_{\l},0)\to (\mathbb
H^{2N+1}_{\l},0)$, $N>n$,
where $L$ denotes the linear embedding as above
and
\begin{equation}\Label{L-}
L_-(z,w):=(z_{\l+1},\ldots, z_n,z_1,\ldots, z_\l,0,-w),
\end{equation}
cannot be transformed into each by composing to the left with
$T\in\Aut(\mathbb H^{2N+1}_{\l},0)$  (cf.\ e.g.\ the proof of
Theorem \ref{t0.3}). However, the following holds.

\bt\Label{t.extra} Let $n$ be even and $H\colon (\mathbb
H^{2n+1}_{n/2},0)\to (\mathbb H^{2(n+k)+1}_{n/2},0)$ a formal
embedding. If $k<n$, then there exists $T\in \Aut(\mathbb
H^{2(n+k)+1}_{n/2},0)$ such that $T\circ H$ is either the linear
embedding $L$ or the embedding $L_-$ given by \eqref{L-}. \et

We would like to point out that the conclusion of
Theorem~\ref{t0.4} fails in general when $k_1+k_2\geq n$. Indeed,
in the case $k_1=0$, $k_2=n$, $\l=0$, the linear embedding of the
sphere can not be represented by a composition of any $T\in
\Aut(\5 H^{2(n+k_2)+1},0)$ and the Whitney embedding; cf.\
\cite{Fa86,Hu1}). It can also fail in the case where the
hyperquadric $\5H^{2N+1}_\l$ is replaced by one with a greater
signature. Such an example is given by the following:

\be\Label{e0.5} Let $k_1=0$, $k_2=2$, $\l=0$, and $M=\5 H^{2n+1}$.
Then, for any formal power series $\phi(z,w)$, the map
$$
H(z,w):=(\phi(z,w),\phi(z,w),z,w)\in \C\times\C\times \C^n\times
\C
$$
embeds $M$ in $\5 H^{2(n+k_2)+1}_{\l'}$ with $\l'=1$. Clearly, not
every such mapping can be represented by a composition of some
$T\in \Aut(\5 H^{2(n+k_2)+1}_{\l'},0)$ and the linear embedding
$(z,w)\mapsto (0,0,z,w)$. \ee

However, the situation in Example~\ref{e0.5} is essentially the only
way in which the conclusion of Theorem~\ref{t0.4} can fail
in case of different signatures.
The reader is referred to Theorem~\ref{t4.2} (and the subsequent remark)
below for the precise statement.

We should point out that if $M\subset\C^{n+1}$ is a smooth
Levi-nondegenerate hypersurface and $H_j\colon M\to \5
H^{2(n+k_j)+1}_\l$, $j=1,2$, are smooth embeddings (rather than
formal), then the conclusion of Theorem~\ref{t0.4} (and,
similarly, of Theorem \ref{t.extra}) is, {\it a priori}, only that
the Taylor series of $T\circ L\circ H_1$ and $H_2$ are equal at
$p$. However, it then follows that $T\circ L\circ H_1$ and $H_2$
are equal also as smooth CR mappings. In the case $\l=0$, this
follows by applying a finite determination result due to the first
author and B. Lamel \cite[Corollary~4]{EL},
%in the case
%$\l=0$ (i.e.\ when $M$ and the hyperquadrics are strictly pseudoconvex)
and in the case $\l>0$, when the Levi form has eigenvalues of both
signs, by using the well-known holomorphic extension of smooth CR
mappings to a neighborhood of $p$.
%in the case $\l>0$.
In particular, Theorem~\ref{t0.4} immediately yields the following corollary.

\bc\Label{cor-rigid}
Let $M_1,M_2$ be compact smooth CR manifolds of hypersurface type
admitting smooth CR embeddings into the unit sphere in $\C^{n+k+1}$ with $k<n/2$.
Then $M_1$ and $M_2$ are globally CR equivalent if and only if
their germs  at some points $p_1\in M_1$ and $p_2\in M_2$ are locally CR equivalent.
\ec

As an application of Corollary~\ref{cor-rigid},
we obtain remarkable families of pairwise locally inequivalent
compact homogeneous hypersurfaces:

\be\Label{lie}
The following CR manifolds naturally appear as $SO(m)$-orbits
contained in the boundaries of Lie balls.
For $m\ge 3$ and $0<\lambda<1$, set
$$M_\lambda:=\{|z_1|^2+\cdots+|z_m|^2 = 1+\lambda^2, \, z_1^2+\cdots+z_m^2=2\lambda\} \subset\C^m.$$
Then $M_\lambda$ are CR submanifolds of $\C^m$ of hypersurface type
and it was shown by W.~Kaup and the third author \cite[Proposition~13.4]{KZ}
that they are pairwise globally CR inequivalent (even pairwise not CR homeomorphic).
% $\lambda_1\ne\lambda_2$, $M_{\lambda_1}$ and $M_{\lambda_2}$
In view of Corollary~\ref{cor-rigid} we conclude that,
if $m\ge 5$, these CR manifolds are even locally CR inequivalent.
\ee

The case of smooth CR embeddings in Theorem~\ref{t0.4}, for $\l=0$
and $k_1=k_2$, was also treated, by a different method, in the
authors' recent paper \cite{EHZ} and previously, for $k_1=k_2=1$,
by Webster \cite{We79}. The case $k_1=0$, i.e.\  of mappings
between hyperquadrics, is of special interest. This situation, for
$\l=0$, has been studied by
 Webster \cite{We79}, Cima-Suffridge \cite{CS83}, Faran \cite{Fa86},
 Forstneri\v c \cite{Fo89},
and the second author \cite{Hu1}. We should mention that in many
of the later papers treating the case $k_1=0$ and $\l=0$ (mappings
between balls), the focus has been on the study of CR mappings
with {\em low initial regularity} (see also the survey
\cite{Hu2}). We shall not address this issue in the present paper.

Theorem \ref{t0.4} can also be viewed as interpolating between the
extreme cases $k_1=k_2\leq n/2$ and $k_1=0$, $k_2<n$ (both of
which had been previously studied in the strictly pseudoconvex
case; see the remarks above). We illustrate this with an example.

\begin{Exa} Let $M\subset\bC^{n+1}$ be a nonspherical ellipsoid,
i.e.\ defined in real coordinates $Z_j=x_j+iy_j$, $j=1,\ldots,
n+1$, by
$$
\sum_{j=1}^{n+1}
\bigg(\frac{x_j}{a_j}\bigg)^2+\bigg(\frac{y_j}{b_j}\bigg)^2=1,
$$
with $a_j\neq b_j$ for, at least, some $j$. After a simple scaling
(if necessary), the ellipsoid $M$ can be described by the equation
$$
\sum_{j=1}^{n+1}\big(A_jZ_j^2+A_j\bar Z_j^2+|Z_j|^2\big)=1,
$$
where $0\leq A_1\leq \ldots\leq A_{n+1}< 1$; the fact that $M$ is
nonspherical means that at least $A_{n+1}>0$. It was observed in
\cite{We78,We99} that the polynomial mapping (of degree
two) $G\colon \bC^{n+1}\to\bC^{n+2}$ given by
$$
G(Z)=\bigg(Z,i\Big(1-2\sum_{j=1}^{n+1} A_jZ_j^2\Big)\bigg)
$$
sends $M$ into the Heisenberg hypersurface $\mathbb
H^{2(n+1)+1}\subset \bC^{n+2}$. Let $p\in M$ and $S$ be any
automorphism of $\mathbb H^{2(n+1)+1}$ which sends $G(p)$ to $0$.
A direct application of Theorem \ref{t0.4} (with $k_1=1$ and
$k_2=k$) then shows that any formal embedding $H\colon (M,p)\to
(\mathbb H^{2(n+k)+1},0)$, with $k<n-1$, must be of the form
$H=T\circ L\circ S\circ G$, where $L$ is the linear embedding
$(\mathbb H^{2(n+1)+1},0)\to(\mathbb H^{2(n+k)+1},0)$ and $T\in
\Aut(\5 H^{2(n+k)+1}_{\l'},0)$. In particular any such formal
embedding is a rational mapping (of degree two).
\end{Exa}

The main technical result in this paper (Theorem~\ref{t1.2} below)
is a statement about the unique solvability of a linear operator
introduced by S. S. Chern and J. K. Moser \cite{CM} in their study
of normal forms for Levi-nondegenerate hypersurfaces. Using this
result, we prove a uniqueness result (Theorem~\ref{t3.2}) for
normalized (formal) mappings between real hypersurfaces given by
equations of the form \eqref{0.3}, where $A$ is actually allowed
to be in a more general class of formal power series than the
class $\3 H_k$ above. Theorem~\ref{t3.2} is the basis for proving
the main results presented above.

The paper is organized as follows. In section~\ref{s1}, we discuss
the Chern-Moser operator, introduce natural classes of formal
power series containing the above classes $\3H_k$ and formulate
our main result Theorem~\ref{t1.2} regarding unique solvability of
this operator for the introduced class. The proof of
Theorem~\ref{t1.2} is given in section~\ref{s2}. A discussion of
the equivalence problem and its explicit solution for the class
$\3H_k$ is given in section~\ref{s3}. The proof of Theorems~
\ref{t0.3} is also given in section~\ref{s3}. In section~\ref{s4},
we study embeddings into hyperquadrics and prove
Proposition~\ref{p0.2}, Theorems~\ref{t0.4} and \ref{t.extra} on
rigidity and the generalization for different signatures in the
form of Theorem~\ref{t4.2}.

\section{Admissible classes and unique solvability of the Chern-Moser operator}
\Label{s1}

In the celebrated paper \cite{CM}, S. S. Chern and J. K. Moser
constructed
%(among other things)
a formal (or convergent in the
case of a real-analytic hypersurface) normal form for Levi-nondegenerate hypersurfaces.
In doing so, they introduced a linear
operator $\3 L$, which we shall refer to as the {\em
Chern-Moser operator}, sending $(n+1)$-tuples $(f,g)=(f_1,\ldots, f_{n},g)$
of $\C^n\times\C$-valued complex formal power
series in $(z,w)\in \bC^{n}\times\bC$
%(or holomorphic functions)
to real formal power series in $(z,\bar z, u)\in\C^n\times\C^n\times\bR$
%, on $\5H^{2n+1}_\l$
%(or real-analytic functions),
defined by
$$
\3L(f,g):= \Im\big(g(z,w) -2i\big\<\bar z,f(z,w) \big\>_{\ell}\big)|_{w=u+i\<z,\bar z\>_\l}\, ,
$$
where
%$\5H^{2n+1}_\l$ and
$\<\cdot,\cdot\>_\l$
is defined by
%\eqref{0.1} and
\eqref{0.2}.
% respectively.
We shall think of the signature
 $\ell$ as being fixed and suppress the dependence of $\3 L$ on
 $\ell$ in the notation.
%Clearly $f(z,u+i\<z,\bar  z\>_\ell)$ is the trace of the function $f$ on the hyperquadric $\5 H^{2n+1}_\ell$
% (given by \eqref{0.1}).
The fundamental importance of the Chern-Moser operator is that a
formal normal form for Levi-nondegenerate hypersurfaces can be
produced (as will be explained below) through the study of the
unique solvability of ${\3 L}$ over certain spaces of formal
power series.

Let us fix $n\geq 1$ and choose coordinates $(z,w)\in
\bC^{n}\times \bC$ with $z=(z_1,\ldots,z_{n})$. We shall also
write $w=u+iv\in\bR \oplus i\bR$.
Let ${\3 O}$ denote the space of real-valued formal power series
in $(z,\bar{z},w,\bar w)$ vanishing of order at least $2$ at $0$.
In particular, any formal power series in $(z,\bar z,u)$
can be viewed as an element in $\3O$ via the substitution $u=(w+\bar w)/2$.
% to weighted order at least $3$ (i.e.\ the lowest nontrivial weighted homogeneous term has
%order at least 3).
Also, let $\3 F$ denote the space of
$(n+1)$-tuples $(f,g)$ of complex formal power series in $(z,w)$
satisfying the following normalization condition
(which characterizes the unit element of the automorphism group of
${\5 H}_{\l}^{2n+1}$ among those of the form $(z,w)\mapsto (z+f(z,w),w+g(z,w)$):
\begin{equation}\Label{1.1}
\frac{\p(f,g)}{\p(z,w)}(0)=0,\quad
\re \Big(\frac{\partial ^2 g}{\partial
w^2}\Big)(0)=0.
\end{equation}
It is not difficult to see that
$\3 L$ sends $\3 F$ into $\3 O$.

By an {\em admissible class}
we shall mean a subset ${\3 S}\subset {\3 O}$
such that the equation
$${\3 L}(f,g)(z,\bar z,u)=\big(A_1(z,\bar z,w,\bar w)-A_2(z,\bar z,w,\bar w)\big)|_{w=u+i\<z,\bar z\>_\l},
\quad (f,g)\in \3F,\quad A_1, A_2\in {\3 S},$$
has only the
trivial solution $(f,g,A_1,A_2) \equiv 0$.
%(Therefore, it must hold that $A_1= A_2$ on $\5H^{2n+1}_\l$).
Admissible classes can be useful in the study of
equivalence problems:
 Let $\3S$ be an admissible class
and $(M_1,0)$ and $ (M_2,0)$ be germs of
(formal) Levi-nondegenerate real hypersurfaces at $0$ defined
respectively by equations of the form
$$v=\<z,z\>_{\l}+A_j(z,\bar{z},w,\bar w),\quad A_j\in\3S,\quad j=1,2.$$
 Then it follows from the proof of Theorem~\ref{t3.2} below that, if there is a
 (formal) holomorphic mapping of the form $H=\id+(f,g)$ with $(f,g)\in \3 F$
 sending $M_1$ to $M_2$,  it must be the identity
map (and hence $M_1=M_2$).
%Moreover, this can clearly only occur when $F_1\equiv F_2$.
We should also point out that if $H$ is any (formal)
biholomorphic mapping (not necessarily normalized as above)
sending $M_1$ to $M_2$, then there exists a unique automorphism
$T\in\Aut(\5 H^{2n+1}_\ell,0)$ such that $H$ factors
as $H=T\circ \2H$, where $\2H$ is normalized as above
(cf.\ e.g.\ Lemma~\ref{l4.1} for a
slightly more general statement). Thus, if there are two mappings
$H$ and $G$ as above, both sending $M_1$ to $M_2$, and which factor through the
same automorphism $T$, i.e.\ $H=T\circ \2H$ and $G= T\circ \2G$
with $\2H$ and $\2G$ normalized as above, then $H$ and $G$ must be equal
(since $H^{-1}\circ G$ will satisfy \eqref{1.1} and send $M_1$ to
itself).  Hence, in this case a defining equation of the form \eqref{0.3}
with $A\in \3 S$
is unique up to an action by the finite-dimensional stability
group $\Aut(\5 H^{2n+1}_\ell,0)$. However, in general, the
action of $\Aut(\5 H^{2n+1}_\ell,0)$ can be complicated, since
an additional renormalization may be required to put $M$
in the form \eqref{1.3} with $A\in \3 S$
after applying an automorphism $T\in \Aut(\5 H^{2n+1}_\l,0)$.

If an admissible class ${\3 S}$ has the property that, for any
$A\in {\3 O}$,
the equation
$${\3 L}(f,g)=A \mod {\3 S}$$
is always solvable with $(f,g)$ satisfying \eqref{1.1},
then ${\3 S}$ is  called a {\em normal space}.
It follows from the definition of an admissible class
that such a solution $(f,g)$ must be unique.
In the case
when ${\3 S}$ is a normal space, any Levi-nondegenerate
hypersurface of
signature $\l$ can be transformed into a
%normal
hypersurface defined by \eqref{1.3} with $A\in {\3 S}$.
In \cite{CM} the authors construct an explicit normal space, which we call the
{\em Chern-Moser normal space} and denote by ${\3 N}$.
Moreover, they prove there that, for a real-analytic $M$,
its normal form and the associated
transformation map are given by convergent power series.

In this paper,
%partially motivated by the author's previous work in [Hu1],
we shall introduce a new admissible class, which is
different from the Chern-Moser normal space, and study the unique
solvability property of the Chern-Moser operator for this class.
%This set is a little more general than needed to prove
%Theorem~\ref{t0.3}, but has the advantage of yielding, without much extra
%work, the more general result given by Theorem~\ref{t3.2}.
Throughout this paper, for a formal power series $A(z,\bar z,w,\bar w)$, we shall use the expansion
\begin{equation}\Label{expand-s}
A(z,\bar z,w,\bar w)=\sum_{\mu,\nu,\gamma,\delta}
A_{\mu\nu\gamma\delta}(z,\bar z)w^\gamma \bar w^\delta,
\quad (z,w)\in\C^n\times\C,
\end{equation}
where $A_{\mu\nu\gamma\delta}(z,\bar z)$
is a bihomogeneous polynomial in $(z,\bar z)$ of bidegree $(\mu,\nu)$
for every $(\mu,\nu,\gamma,\delta)$.

\bd\Label{d1.1} For every positive integer $k$, define $\3
S_{k}\subset\3O$ to be the subset of all real-valued formal power
series $A(z,\bar z,w,\bar w)$ with $(z,w)\in \C^n\times\C$,
satisfying either of the following equivalent conditions:
\begin{enumerate}
\item[(i)] $A_{\mu\nu\gamma\delta}=0$ if $\nu+\delta\le 1$ and $\R(A_{\mu\nu\gamma\delta})\le k$ if $\delta\le 1$;
\item[(ii)] for each fixed $(\mu,\nu,\gamma,\delta)$ with $\delta\le 1$,
we can write
$$A_{\mu\nu\gamma\delta}(z,\bar z)w^\gamma\bar w^\delta=\sum_{j=1}^k \phi_j(z,w)\1{\psi_j(z,w)}$$
for some holomorphic polynomials $\phi_j$, $\psi_j$ with no constant or linear terms.
\end{enumerate}
\ed

Note that, since $A$ is real-valued, condition (i) implies also
$A_{\mu\nu\gamma\delta}=0$ if $\mu+\gamma\le 1$ and
$\R(A_{\mu\nu\gamma\delta})\le k$ if $\gamma\le 1$.
Our main result regarding the set $\3S_k$ is the following
uniqueness property.

\bt
\Label{t1.2}
Given any $A\in \3S_{n-1}$, the equation
\begin{equation}\Label{lfg}
{\3 L}(f,g)(z,\bar z,u)=A(z,\bar z,w,\bar w)|_{w=u+i\<z,\bar z\>_\l},
\end{equation}
for $(f,g)$ satisfying \eqref{1.1}
has only the trivial solution $(f,g)\equiv 0$.
\et

The conclusion of Theorem~\ref{t1.2} reduces to a fundamental
result of Chern-Moser \cite{CM}, when the right-hand side in
\eqref{lfg} is in their normal space  ${\3 N}$. Here, we recall
that ${\3 N}$ consists of all formal power series
$\sum_{\a,\b}A_{\a\-{\b}}(u)z^{\a}\bar{z}^{\b}$ with $|\a|,
|\b|\ge 2$, and $A_{2\-{2}}, A_{2\-{3}}, A_{3\-{3}}$ satisfying
certain trace conditions as described in \cite[pp 233]{CM}.
%, cf.\ also \{} [pp 315, BER2]).
It is clear that the right-hand side of \eqref{lfg}, with $A\in
\mathcal S_{n-1}$, is in general not in the Chern-Moser normal
space. For instance, $A(z,\bar z,w,\bar w): =|w|^4$ is easily seen
to be in $\3 S_1$ whereas $A|_{w=u+\<z,\bar z\>_\l}$ is not in the
Chern-Moser normal space $\3N$.

\section{Proof of Theorem~\ref{t1.2}}\Label{s2}
The proof of Theorem~\ref{t1.2} is based on the following lemma from
\cite{Hu1}.

\bl\Label{l2.1} Let $\phi_j$, $\psi_j$, $j=1,\ldots, n-1$, be
germs at $0$ of holomorphic functions in ${\5 C}^n$ and
$H(z,\bar{z})$ a germ at $0$ of a real-analytic function
satisfying
\begin{equation}\Label{2.1.1}
H(z,\bar \xi)\<z,\bar \xi\>_{\l} =
\sum_{j=1}^{n-1}\phi_j(z)\-{\psi_j(\xi)}.
\end{equation}
Then
$$H(z,\bar \xi) =
\sum_{j=1}^{k}\phi_j(z)\-{\psi_j(\xi)} = 0.$$
\el

Lemma~\ref{l2.1} is the content of Lemma~3.2 in \cite{Hu1} when $\l=0$.
Also,
in the statement of \cite[Lemma 3.2]{Hu1}, it is assumed that
$\phi_j(0)=\psi_j(0)=0$. However, the argument there can be used
for the proof of Lemma~\ref{l2.1} without any change. We should also
point out that in Lemma~\ref{l2.1} it is enough to assume the identity
$$H(z,\bar z)\<z,\bar z\>_{\l}=
\sum_{j=1}^{k}\phi_j(z)\-{\psi_j(z)},$$
since the
identity \eqref{2.1.1} then follows by a standard complexification
argument. Using Lemma~\ref{l2.1} and a simple induction argument, which we
leave to the reader, we obtain the following generalization of
Lemma~\ref{l2.1}.

\bl\Label{l2.2} Let $\phi_{jp},\ \psi_{jp}$, $j=1,\ldots, n-1$,
$p=0,\ldots,q$, be germs at $0$ of holomorphic functions in
${\5C}^n$ and $H(z,\bar{z})$ a germ at $0$ of a real-analytic
function satisfying
\begin{equation*}%\Label{2.2.1}
H(z,\bar \xi)
\<z,\bar \xi\>_{\l}^{q+1}
=
\sum_{p=0}^q
\Big(
\sum_{j=1}^{n-1}
\phi_{jp}(z)\-{\psi_{jp}(\xi)}
\Big)
\<z,\bar\xi\>_\ell^p.
\end{equation*}
Then
$$H(z,\bar \xi) = \sum_{j=1}^{n-1}\phi_{jp}(z)\-{\psi_{jp}(\xi)} = 0, \quad p=1,\ldots, q.$$
\el

\bpf[Proof of Theorem $\ref{t1.2}$] Recall that the Segre variety
$Q_{(\1\xi,\1\eta)}$ of ${\5 H}_{\l}^{2n+1}$ with respect to
$(\1\xi,\1\eta)\in {\5C}^{n}\times {\5 C}$ is the complex
hyperplane defined by
$$Q_{(\1\xi,\1\eta)}:=\{(z,w)\in {\5 C}^{n}\times {\5 C}:
w-{\eta}=2i\<z,{\xi}\>_{\l}\}.$$
If we set
\begin{equation}\Label{vf}
L_j=\frac {\p}{\p z_j}+2i\d_j {\xi_j}\frac{\p}{\p w},
\end{equation}
 with $\d_j =-1$ for $j\le
\l$ and $\d_j=1$ for $j>\ell$, then $(L_j)_{1\le j\le n}$ forms a
basis of the space of $(1,0)$-vector fields along
$Q_{(\1\xi,\1\eta)}$. By noticing that the equation \eqref{lfg}
%${\3 L}(f,g)=S|_{\5H^{2n+1}_\l}$
is the same as
$$g(z,w)-\-{g(z,w)}-2i\<\bar{z},f(z,w)\>_{\l}-2i\<z,\-{f(z,w)}\>_{\l}
=2i A(z,\bar z,w,\bar w), \quad (z,w)\in \5 H^{2n+1}_\ell,$$
and then complexifying in the standard way,
we obtain
\begin{equation}\Label{2.1}
g(z,w)-\1g(\xi,\eta)-2i\<\xi,f(z,w)\>_\l
-2i\<z,\1f(\xi,\eta)\>_\l
=2i A(z,\xi,w,\eta)
\end{equation}
which holds for $(z,\xi,w,\eta)$ satisfying
$w-{\eta}=2i\<z,{\xi}\>_{\l}$ (or, if we think of $(\xi,\eta)$
as being fixed, for $(z,w)\in Q_{(\1\xi,\1\eta)}$).
We shall also use the identity obtained
by applying $L_j$ to \eqref{2.1}:
\begin{equation}\Label{2.2}
L_j(g(z,w))-2i\<{\xi},L_j(f(z,w))\>_{\l}
-2i\1f_j(\xi,\eta)= 2i L_j A(z,\xi,w,\eta),
\quad w={\eta}+2i\<z,{\xi}\>_{\l}.
\end{equation}

In view of \eqref{1.1} we have the expansions
\begin{equation}\Label{expand}
f(z,w)=\sum_{\mu+\nu\ge 2} f_{\mu\nu} (z) w^\nu, \quad
g(z,w)=\sum_{\mu+\nu\ge 2} g_{\mu\nu} (z) w^\nu,
%\quad
%S(z,\xi,w,\eta)=\sum_{{\mu+\gamma\ge 2\atop \nu+\delta\ge 2}}
%S_{\mu\nu\gamma\delta}(z,\xi)w^\gamma\eta^\delta
\end{equation}
where $f_{\mu\nu} (z)$ and $g_{\mu\nu} (z)$
are homogeneous polynomials of degree $\mu$.
We also use the expansion \eqref{expand-s}.
We allow the indices to be arbitrary integers
by using the convention that all coefficients not appearing in
\eqref{expand-s} and \eqref{expand} are zero. We now let $w=0$,
$\eta=\eta(z,\xi)=-2i\<{z},\xi\>_{\l}$ in \eqref{2.1} and \eqref{2.2},
and equate bihomogeneous terms in $(z,\xi)$ of a fixed bidegree $(\a,\b)$.
Comparing terms with $\b=1$ and $\b=0$ in \eqref{2.1} we obtain
\begin{equation}\Label{1st}
f(z,0)=0, \quad g(z,0)=0.
\end{equation}
For terms of a fixed bidegree $(\a,\b)$ with $\b\ge 2$, we have
\begin{equation}\Label{2nd}
-\1g_{\b-\a,\a}(\xi)\eta^\a
-2i \<z,\1f_{\b-\a+1,\a-1}(\xi)\eta^{\a-1} \>_\l
=2i \sum_{p=0}^{\a-2} A_{\a-p,\b-p,0,p}(z,\xi)\eta^p,
\end{equation}
where $\eta=-2i\<{z},\xi\>_{\l}$.
Since $\R(A_{\mu \gamma 0\delta})<n$ for all $(\mu,\gamma,\delta)$,
Lemma~\ref{l2.2} implies $A_{\mu \gamma 0\delta}= 0$
for $\eta=-2i\<{z},\xi\>_{\l}$ and hence
\begin{equation}\Label{2nd'}
\1g_{\b-\a,\a}(\xi)\eta+
2i \<z,\1f_{\b-\a+1,\a-1}(\xi) \>_\l
=0, \quad \eta=-2i\<{z},\xi\>_{\l}, \quad \b\ge 2.
\end{equation}
In particular, in view of \eqref{vf} and \eqref{1st}, we have
\begin{equation}\Label{ders}
(L_j(f,g)) (z,0) = 2i\delta_j \xi_j \sum_\mu (f_{\mu 1},g_{\mu 1})(z,0), \quad
(L_j A)(z,\xi,0,\eta) =
2i\delta_j \xi_j
\sum_\mu A_{\mu \gamma 1\delta} (z,\xi)\eta^\delta.
\end{equation}

We now apply the same procedure (i.e.\
collect terms of a fixed bidegree $(\a,\b)$ in $(z,\xi)$) to \eqref{2.2}
using \eqref{ders}.
For $\b=0$ we then obtain no terms and for $\b=1$ and $2$ respectively the
identities
\begin{equation}\Label{b0}
g_{\a 1}(z)= 0,
\quad
\<\xi, 2i\delta_j \xi_j f_{\a 1}(z) \>_\l
+ \1f_{j;2-\a,\a}(\xi)\eta^\a =0,
\quad \eta=-2i\<{z},\xi\>_{\l},
\end{equation}
where we use the notation $f_{\mu\nu}=f_{\mu,\nu}=(f_{1;\mu,\nu},\ldots,
f_{n;\mu,\nu})$.
Finally, for $\b\ge 3$, the same comparison yields
\begin{equation*}%\Label{}
-2i \1f_{j;\b-\a,\a}(\xi)\eta^\a = 2i\delta_j \xi_j \sum_{p=0}^{\a-1}
A_{\a-p,\b-p-1,1,p}(z,\xi)\eta^p.
\end{equation*}
Using the assumption $\R(A_{\mu \gamma 1\delta})<n$
and applying Lemma~\ref{l2.2} as above, we conclude
that $f_{\mu\nu}(z)=0$ for $\mu+\nu\ge 3$.
Substituting this into \eqref{2nd'} we conclude that $g_{\mu\nu}(z)=0$ for
$\mu+\nu\ge 3$.
Also, substituting $f_{21}(z)=0$ into \eqref{b0} yields $f_{02}=0$.
Since $f_{20}=0$ in view of \eqref{1st},
the identity \eqref{2nd'} for $(\a,\b)=(1,2)$ implies $g_{11}=0$.

It remains to show that $(f_{11}(z),g_{02})=0$,
where we drop the argument $z$ for $g_{02}$ since the latter is a constant.
Rewriting \eqref{2nd'} for $\a=\b=2$ and the second identity
in \eqref{b0} for $\a=1$, we obtain respectively
\begin{equation}\Label{g02}
\1g_{02}\<z,\xi\>_\l=\<z,\1f_{11}(\xi) \>_\l,
\quad
\delta_j \xi_j \<\xi,f_{11}(z) \>_\l = \1f_{j;11}(\xi) \<z,\xi \>_\l.
\end{equation}
Setting $\xi=\bar z$ and using the normalization \eqref{1.1},
\eqref{g02} implies
\begin{equation}\Label{11id}
\Re\<z,\1{f_{11}(z)}\>_\l=0, \quad
\delta_j \bar z_j \<\bar z,f_{11}(z) \>_\l = \1{f_{j;11}(z)} \<z,\bar z\>_\l.
\end{equation}
Recall that $f_{j;11}(z)$ is a linear function in $z$
that we write as $f_{j;11}(z)=\sum_k f_j^k z_k$.
Then the second identity in \eqref{11id} can be rewritten as
\begin{equation*}%\Label{}
\delta_j \bar z_j \sum_{s,k} \delta_s \bar z_s f_s^k z_k =
\sum_{l,m} \delta_m \1{f_j^l} \bar z_l z_m \bar z_m,
\end{equation*}
from which we conclude that $f_{j;11}(z)=c_jz_j$ for some $c_j\in\R$,
$j=1,\ldots,n$.
Substituting this into the first identity in \eqref{11id}
we see that $\Re c_j=0$ for all $j$ and hence $f_{11}(z)=0$.
Now the first identity in \eqref{g02} yields $g_{02}=0$.
The proof is complete.
\epf

We mention that in Definition~\ref{d1.1}, it is important to
assume that $\phi_j, \psi_j$ have no linear terms in $(z,w)$.
Otherwise, the conclusion in Theorem~\ref{t1.2} fails as can be
easily seen by the example  $f(z,w):=(0,\ldots,0,-\chi(z))$,
$g:=0$ and $A:=z_n\1{\chi(z)}+\chi(z)\1{z}_n$ with $\chi(z)$ being
any holomorphic function in $z$ vanishing at $0$. Also the class
$\3S_{n-1}$ in Theorem~\ref{t1.2} cannot be replaced by any
$\3S_k$ with $k\ge n$ as the following example shows.

\be\Label{e2.2}
Set $f:=(0,\ldots,0,\frac {i}{2}z_n w)$, $g:=0$ and $A:=|z_n|^2\<z,\bar z\>_\l$.
Then $(f,g)\ne 0$ satisfies \eqref{1.1} and solves the equation \eqref{lfg}.
\ee

\section{Application to the equivalence problem}\Label{s3}
In the situation of Theorem~\ref{t1.2} we have shown that $(f,g)$
must vanish. In view of \eqref{lfg}, it follows that the
restriction of $A$ to $\5H^{2n+1}_\l$ also vanishes. However, it
is easy to see that the full power series $A(z,\bar z, w,\bar w)$
need not necessarily vanish. In this section we shall refine the
class $\3S_{n-1}$ to a smaller one $\2{\3S}_{n-1}\subset
\3S_{n-1}\subset\3O$ with the property that any $A\in
\2{\3S}_{n-1}$ which vanishes on $\5H^{2n+1}_\l$ vanishes
identically.

\bd\Label{d3.1} For every positive integer $k$, define $\2{\3
S}_{k}\subset\3O$ to be the subset of all real-valued formal power
series $A(z,\bar z,w,\bar w)$ with $(z,w)\in \C^n\times\C$,
satisfying either of the following equivalent conditions:
\begin{enumerate}
\item[(i)] $A_{\mu\nu\gamma\delta}=0$ if $\nu+\delta\le 1$ and $\R(A_{\mu\nu\gamma\delta})\le k$ otherwise;
\item[(ii)] for each fixed $(\mu,\nu,\gamma,\delta)$, we can write
$$A_{\mu\nu\gamma\delta}(z,\bar z)w^\gamma\bar w^\delta=\sum_{j=1}^k \phi_j(z,w)\1{\psi_j(z,w)}$$
for some holomorphic polynomials $\phi_j$, $\psi_j$ with no constant or linear terms.
\end{enumerate}
\ed

The property of $\2{\3S}_{n-1}$ mentioned above is a consequence
of the following statement.

\bl\Label{l3.3}
Let $A(z,\bar{z},w,\bar{w})$ with $(z,w)\in \bC^n\times\bC$
 be a formal power series in the class $\2{\3 S}_{n-1}$.
Assume that
$A(z,\bar{z},w,\bar w)|_{w=u+i\<z,\bar{z}\>_{\l}} \equiv 0$
as a formal power series in $(z,\bar{z},u)$. Then
$A(z,\bar{z},w,\bar{w})\equiv 0$ as a formal power series in
$(z,\bar{z},w,\bar{w})$.
\el

\bpf
We use the expansion \eqref{expand-s} for $A$.
The same complexification argument as in the proof of Theorem~\ref{t1.2} yields
$$
A(z,\xi,w,\eta) =0  \text{ for } w=\eta+2i\<z,\xi\>_\ell
$$
which can be rewritten as
\begin{equation*}%\Label{3.6}
\sum_{\mu,\nu,\gamma,\delta}
A_{\mu\nu\gamma\delta}(z,\xi)(\eta+2i\<z,\xi\>_\ell)^\gamma \eta^\delta\equiv 0.
\end{equation*}
Assume, in order to reach a contradiction, that $A(z,\bar z,w,\bar w)\not\equiv 0$.
Then, there is a smallest nonnegative integer
$\delta_0$ such that $A_{\mu\nu\gamma\delta_0}(z,\xi)\not\equiv 0$ for some
$(\mu,\nu,\gamma)$. By factoring out $\eta^{\delta_0}$ (of course, if
$\delta_0=0$, then we do not need to factor anything)
and setting $\eta=0$, we obtain
\begin{equation*}%\Label{3.7}
\sum_{\mu,\nu,\gamma}
A_{\mu\nu\gamma\delta_0}(z,\xi)\big(2i\<z,\xi\>_\ell\big)^\gamma \equiv 0.
\end{equation*}
Isolating terms of a fixed bidegree $(\a,\b)$ in $(z,\xi)$, we deduce
\begin{equation}\Label{fin-id}
\sum_p A_{\a-p,\b-p,p,\delta_0}(z,\xi)
\big(2i\<z,\bar \xi\>_\ell\big)^p\equiv 0.
\end{equation}
By the definition of the class $\2{\3S}_{n-1}$,
we have $\R(A_{\mu\nu\gamma\delta_0})<n$
for every $(\mu,\nu,\gamma)$.
Hence Lemma~\ref{l2.2}, applied to the identities \eqref{fin-id}
for all $(\a,\b)$, yields $A_{\mu\nu\gamma\delta_0}\equiv 0$
for all $(\mu,\nu,\gamma)$
in contradiction with the choice of $\delta_0$.
This completes the proof of the lemma.
\epf

Note that for $\alpha_1,\alpha_2\in \bR$,
$A_1\in\3 S_{k_1}$,
$A_2\in\3 S_{k_2}$, it holds that
$\alpha_1 A_1+\alpha_2 A_2\in \3 S_{k_1+k_2}$.
Hence Theorem~\ref{t1.2} together with Lemma~\ref{l3.3}
imply that $\3 S_{k}$ is an admissible
class in the sense of section \ref{s1} for $k<n/2$.

Given a formal power series
$A(z,\bar z,w,\bar w)$, we associate to it a formal power series
$A^0(z,\bar{z},u)$, with $u\in\bR$, defined by
$$A^0(z,\bar{z},u):=A(z,\bar{z},w,\bar w)|_{w=u+i\<z,\bar{z}\>_{\l}};$$
in other words, $A^0$ can be viewed as the trace of
$A$ along the hyperquadric $\5 H^{2n+1}_\ell$.
Thus, given a formal power series $A\in \3 S_k$, we can
associate to it two formal hypersurfaces $M^0$ and $M$ in $\C^{n+1}$
 as follows:
\begin{equation}\Label{3.3}
M^0:=\{\im w=\<z,\bar{z}\>_{\l}+A^0(z,\bar{z},u)\}, \quad
 M:=\{\im w=\<z,\bar{z}\>_{\l}+A(z,\bar{z},w,\-{w})\}.
\end{equation}
(Observe that the equation defining $M$  is not in graph form.)
For these classes of hypersurfaces,
Theorem~\ref{t1.2} and Lemma~\ref{l3.3}
can be used for the study of the equivalence problem as follows.

\bt\Label{t3.2}
Let $A_1\in \2{\3 S}_{{k_1}}$, $A_2\in \2{\3 S}_{{k_2}}$, and
 define the formal hypersurfaces $M_j^0$ and $M_j$ in $\bC^{n+1}$ for $j=1,2$,
 by  \eqref{3.3} $($with $A_j$ in the place of $A)$.
Let $$H(z,w)=(z+f(z,w),w+g(z,w))$$ be a formal biholomorphic map
with $(f,g)$ satisfying the normalization condition \eqref{1.1}
and which sends $M_1$ into $M_2$ $($or $M^0_1$ into $M^0_2)$.
%or $M_1$ into $M^0_2
Then, if $k_1+k_2<n$, it must hold that $H(z,w)\equiv (z,w)$ and $A_1\equiv  A_2$.
\et

\bpf
We shall prove the theorem when $H$
maps $M_1$ to $M_2$; the proof of the other case is
similar and left to the reader.
By the implicit function theorem, $M_1$ can
be represented by the equation
$$\im w=\<z,\bar{z}\>_{\l}+\wt{A}_1(z,\bar{z},u),$$
 where $\wt{A}_1$ is uniquely obtained  by solving for $v$ in the
 equation
$$v=\<z,\bar{z}\>_{\l}+A_1(z,\bar{z},u+iv,u-iv).$$
Since the formal map $H=\id+(f,g)=:(F,G)$
sends $M_1$ to $M_2$, we have the following identity
\begin{equation}\Label{3.8}
\im G-\<F,\bar F\>_{\l}=A_2(F,\bar{F},{G},\bar{G}),
%\wt{F}_2^0(f,\-{f},\re (g)),
\end{equation}
when we set
$$
w=u+i\<z,\bar z\>_\ell + iA_1(z,\bar z,w,\bar w)
=u+i\<z,\bar z\>_\ell + i\wt{A}_1(z,\bar z,u).
$$

As in \cite{CM}, we assign the variable $z$ the weight
$1$ and $w$ (as well as $u$ and $v$) the weight $2$. Recall that a
holomorphic polynomial $h(z,w)$ is then said to
be weighted homogeneous of weighted
degree $\sigma$ if $h(tz,t^2w)=t^{\sigma}h(z,w)$ for any complex
number $t$. A real polynomial $h(z,\bar{z},u)$ (resp.\ $h(z,\bar z,w,\bar w)$)
 is said to be weighted homogeneous of weighted degree $\sigma$, if
 for any real number $t$, $h(tz,t\bar{z},t^2u)=t^{\sigma}h(z,\bar{z},u)$
(resp.\ $h(tz,t\bar{z},t^2w,t^2\bar w)=t^{\sigma}h(z,\bar{z},w,\bar w)$).
 In all cases, we  write $\h{deg}_{wt}(h)=\sigma$.
(By convention, $0$ is a weighted
 homogeneous holomorphic
polynomial
 of any degree.)
We shall denote by $h^{(\sigma)}$
the term of weighted degree $\sigma$ in the expansion of $h$
into weighted homogeneous terms.

It is sufficient to prove the following claim:
$(f^{(\tau-1)},g^{(\tau)})\equiv 0$ and $A_1^{(\tau)}\equiv
A_2^{(\tau)}$ for $\tau\ge 1$. For $\tau=1$, this follows directly
from the assumptions. We shall prove the claim for a general
$\tau\geq 1$ by induction. Assume that it holds for $\tau<\sigma$.
Then, using the fact that, on $M_1$,
$$w=u+i\<z,\bar z\>+iA_1
\big(z,\bar z,u+i\<z,\bar z\>+i\wt{A}_1(z,\bar z,u),
u-i\<z,\bar z\> -i\wt{A}_1(z,\bar z,u)\big),$$
and collecting terms of weighted degree $\sigma\ge 2$ in \eqref{3.8},
we obtain
\begin{multline}
\3 L(f^{(\sigma-1)},g^{(\sigma)})(z,\bar z,u)
=\big(
A_2(z,\bar z,u+i\<z,\bar z\>+i\wt{A}_1(z,\bar z,u),u-i\<z,\bar z\>-i\wt{A}_1(z,\bar z,u)
\big)^{(\sigma)} \\
-\big(
A_1(z,\bar z,u+i\<z,\bar z\>+i\2A_1(z,\bar z,u),
u-i\<z,\bar z\>-i\2A_1(z,\bar z,u)
\big)^{(\sigma)}\\
=\big(A_2^{(\sigma)}(z,\bar z,w,\bar w)-A_1^{(\sigma)}(z,\bar z,w,\bar w)\big)|_{w=u+i\<z,\bar z\>_\l}.
\end{multline}
Since $A_2-A_1\in \3 S_{n-1}$ by the assumption, we conclude by
Theorem~\ref{t1.2} that $(f^{(\sigma-1)},g^{(\sigma)})\equiv 0$.
Furthermore, since also $A_2-A_1\in \2{\3 S}_{n-1}$,
Lemma~\ref{l3.3} implies that
 $A_1^{(\sigma)}\equiv A_2^{(\sigma)}$, which completes the induction.
The proof of the theorem is complete. \epf

\bpf[Proof of Theorem $\ref{t0.3}$] Let $M_1$ and $M_2$ be as in
Theorem~\ref{t0.3} and let $H=(F,G)$ be a formal invertible
mapping $(\bC^{n+1},0)\to (\bC^{n+1},0)$ sending $M_1$ to $M_2$;
that is, \eqref{3.8} holds when $(z,w)\in M_1$. By collecting
terms of weighted degree $1$ and $2$, as in the proof of
Theorem~\ref{t3.2}, we see that there are $\lambda>0$, $r\in \bR$,
$a\in \bC^n$, $\sigma=\pm1$ and an invertible $n\times n$ matrix
$U$ such that
$$
\frac{\p (F,G)}{\p (z,w)}(0):=\begin{pmatrix}
F_z(0)& G_z(0) \\
F_w(0) & G_w(0)
\end{pmatrix}=
\begin{pmatrix}
\lambda U & 0 \\
\lambda a U & \sigma \lambda ^2
\end{pmatrix},
\quad
 {\Re}\Big(\frac{\p^2 G}{\p w^2}\Big)(0)=2\sigma\lambda ^2r,
$$
where the linear transformation $z\to z U$ preserves the Hermitian
product $\<\cdot,\cdot\>_\l$ modulo $\sigma$, i.e.\ $\<zU,\bar z
\bar U\>_\l=\sigma\<z,\bar z\>_\l$. Observe that $\sigma$ must be
$+1$ unless $\l=n/2$. Now, let us define $T\in \Aut(\5
H^{2n+1}_\l,0)$ by \eqref{autos0} (See e.g.\ \cite{CM}
% or [BER]
for the fact that $T\in\Aut(\5 H^{2n+1}_\l,0)$.) It is easy to
verify that $\wt{H}:=T^{-1}\circ H$ satisfies the normalization
condition (1.1). Also, a straightforward calculation, which is
left to the reader, shows that $\wt{H}$ maps $M_1$ (via $M_2$)
into a real (formal) hypersurface in $\C^{n+1}$ which can be
defined by
$$\im w=\<z,\bar z\>_{\l}+\2 A_2(z,\bar z,w,\bar w)$$
with
\begin{equation}\Label{3.10}
\2A_2(z,\bar z,w,\bar w) \equiv
\sigma\lambda^{-2} |q(z,w)|^2
\, A_2\circ (T(z,w),\1{T(z,w)}),
\end{equation}
where $q(z,w)$ denotes the denominator in \eqref{autos0}. Observe
that if $A_2$ belongs to $\mathcal H_k$, then, in view of
Definition~\ref{d0.1}, so does $\tilde A_2$. Thus, by our
assumptions, we have $A_1\in \3H_{k_1}$ and also $\2A_2\in
\3H_{k_2}$. Since $\3 H_k\subset\2{\3 S}_k$ for all $k$ as is easy
to see, it follows from Theorem~\ref{t3.2} that $\2H\equiv \id$
and $A_1\equiv \2A_2$. The conclusion of Theorem~\ref{t0.3} now
follows from \eqref{3.10} and the proof is complete. \epf

\section{Embeddings of real hypersurfaces in
hyperquadrics}\Label{s4}

Let $M$ be as in Proposition~\ref{p0.2} and $H\: (\bC^{n+1},p)\to
(\bC^{N+1},0)$ be a formal holomorphic embedding sending $M$ into
$\5 H_{\l'}^{2N+1}$ (with $N\geq n$ and $\l'\leq N/2$). We shall
assume that $H$ is CR transversal (i.e.\  $\partial H_{N+1}(p)\neq
0$, see section \ref{s0}; also, recall that this is automatic when
$\l'<n$). Let us write $H=(F,G)\in\C^N\times\C$. The fact that $H$
sends $(M,0)$ into $(\5 H^{2N+1}_{\l'},0)$ means that
\begin{equation}\Label{embed}
\im G(z,w)=\<F(z,w),\1{F(z,w)}\>_{\l'} \text{ for } (z,w)\in M.
\end{equation}
By collecting terms of weighted degree 1 and $2$ as in the proofs
of Theorem~\ref{t0.3} and \ref{t3.2} and using the CR
transversality assumption, we see that
\begin{equation}\Label{H-jac}
\frac{\p (F,G)}{\p (z,w)}(0)=
\begin{pmatrix}
\lambda V & 0 \\
\lambda a  & \sigma \lambda ^2
\end{pmatrix},
\end{equation}
for some $\sigma=\pm 1$, $ \lambda>0$, $a\in\C^N$, and an $n\times
N$ matrix $V$ of rank $n$ satisfying
\begin{equation}\Label{productid}
\<zV,\bar z\bar V\>_{\l'} = \sigma \<z,\bar z\>_\l.\end{equation}
Observe that, in view of \eqref{productid}, if $\sigma=-1$ then
$\l' \geq n-\l$. If $\l'<n-l$, then we must have $\sigma=1$. It
will be convenient to renumber the coordinates $z'=(z'_1,\ldots,
z'_N)\in\C^{N}$, according to the sign of $\sigma$, so that
\begin{equation}\Label{4.1}
\<z',\bar z'\>_{\l'}= -\sigma \sum_{j=1}^{\l}|z'_j|^2+\sigma
\sum_{j=\l+1}^{n}|z'_j|^2
-\sum_{j=n+1}^{n+s}|z'_j|^2+\sum_{j=n+s+1}^{N}|z'_j|^2,
\end{equation}
where $s=\l'-\l$ if $\sigma=1$ and $s=\l'-(n-\l)$ if $\sigma =-1$.
We choose coordinates $(z,w)\in \bC^n\times\bC$, vanishing at $p$,
such that $M$ is defined by an equation of the form \eqref{0.3},
where $A$ vanishes to fourth order at $0$. We shall need the
following lemma.

\bl \Label{l4.1} Let $M\subset\bC^{n+1}$ be defined by an equation
of the form \eqref{0.3} and $H=(F,G)\: (\bC^{n+1},p)\to
(\bC^{N+1},0)$ a formal holomorphic, CR transversal embedding
sending $(M,0)$ into $\5 H_{\l'}^{2N+1}$, with $N \ge n$. Let
$\sigma=\pm 1$ so that \begin{equation}\Label{sigma}\sigma
\frac{\partial G}{\partial w}(0)>0\end{equation} and renumber the
coordinates $z'\in \bC^N$ such that $\<z',\bar z'\>_{\l'}$ is
given by \eqref{4.1}. Then, there exists $T\in\Aut(\5
H^{2N+1}_{\l'},0)$ such that $T^{-1}\circ H$ is of the form
\begin{equation}\Label{form}
(z,w)\mapsto (z+f(z,w),\phi(z,w),\sigma w+g(z,w))\in
\C^n\times\C^{N-n}\times\C
\end{equation}
with $d\phi(0)=0$ and  $(f,g)$ satisfying the normalization conditions
\eqref{1.1}.
\el

\begin{Rem} Recall that if $\l'<n-\l$ then, as observed above, we
must have $\sigma=+1$.
\end{Rem}

\bpf As mentioned above, $H=(F,G)$ satisfies \eqref{H-jac}, where
$\lambda>0$, $a\in\bC^N$, and $V$ satisfies \eqref{productid}. Let
$r\in \bR$ be such that
$$
 {\Re}\Big(\frac{\p^2 G}{\p w^2}\Big)(0)=2\lambda ^2r.
$$
By standard linear algebra, we can extend $V$ to an $N\times N$
matrix $U$ such that its first $n$ rows are the those of $V$ and
such that the linear transformation $z'\mapsto z' U$ preserves the
Hermitian form $\<\cdot,\cdot\>_{\l'}$, i.e.\ $ \<z'U,\bar z'\bar
U\>_{\l'}=\<z',\bar z'\>_{\l'}$. (Recall that
$\<\cdot,\cdot\>_{\l'}$ is given by \eqref{4.1}.) Now define
$T\:(\bC^{N+1},0)\to (\bC^{N+1},0)$, as in the proof of
Theorem~\ref{t0.3} above, to be the element of $\Aut(\5
H^{2N+1}_{\l'},0)$ given by
\begin{equation}\Label{autos1}
T(z',w'):= \frac {(\lambda(z'+bw')U,\sigma \lambda^{2}w')}
{1-2i\<z',\bar {b}\>_{\l'} -(r+i\<b,\bar{b}\>_{\l'})w'},
\end{equation}
where $b=\sigma a U^{-1}$. A straightforward calculation, which is
left to the reader, shows that $T$ satisfies the conclusion of the
lemma. \epf

\bpf[Proof of Proposition $\ref{p0.2}$] In view of Lemma
\ref{l4.1}, we may assume, after composing $H$ with some $T\in
\Aut(\5 H^{2N+1}_{\l'},0)$, that $H$ is of the form \eqref{form}
with $d\phi(p)=0$ and $(f,g)$ satisfying \eqref{1.1}. Then the
mapping $H^0(z,w):=(z,\sigma w) + (f(z,w),g(z,w))$ is a local
biholomorphism $(\bC^{n+1},p)\to (\bC^{n+1},0)$ and therefore we
may consider a formal change of coordinates given by $(\hat
z,\sigma \hat w):= H^0(z,w)$ or, equivalently, $(\hat z,\hat
w):=(z+f,w+\sigma g)$; by a slight abuse of notation, we shall
drop the $\hat{}$ over the new coordinates $(\hat z,\hat w)$ and
denote also these coordinates by $(z,w)$. The fact that $H$ is an
embedding sending $M$ into $\5 H^{2N+1}_{\l'}$ means, in view of
\eqref{embed}, that, in the new coordinates, $M$ can be defined by
the equation \eqref{0.3}, where
\begin{equation}\Label{4.4}
A(z,\bar z,w,\bar w):= -\sigma \sum_{j=1}^{s} |\tilde
\phi_j(z,w)|^2 +\sigma \sum_{j=s+1}^{N-n} |\tilde \phi_j(z,w)|^2,
\quad \tilde \phi_j:=\phi_j\circ (H^0)^{-1},
\end{equation}
where $s=\l'-\l$ if $\sigma=1$ and $s=\l'-(n-\l)$ if $\sigma=-1$.
(Cf. also \cite{Fo86}.) By Definition~\ref{d0.1}, the fact that
$A$ is given by \eqref{4.4} means that $\R(A)\leq N-n$ (but not
necessarily that $\R(A)=N-n$) and hence $A\in \3 H_{N-n}$. Also,
observe that the linear subspace $V_A$ defined in the introduction
is an $\R(A)$-dimensional subspace of the span, over $\bC$, of the
formal series $\tilde\phi_1,\ldots,\tilde\phi_{N-n}$. By choosing
a basis $\psi_1,\ldots,\psi_{\R(A)}$ for $V_A$ such that $A$ is
given by the analogue of \eqref{0.4} with $\psi_j$ instead of
$\tilde \phi_j$ (cf.\ section~\ref{s0}), one can conclude using
standard linear algebra that $\s(A)\leq \l'-\l$ when $\sigma=+1$
and $\s(A)\leq  N-\l'-\l$ when $\sigma=-1$. Since $\sigma=-1$ is
only possible if $\l\geq n-l$, this proves one implication in
Proposition~\ref{p0.2}.

For the other implication, we must show that given coordinates
$(z,w)\in \bC^n\times \bC$ vanishing at $p$ such that $M$ is
defined by \eqref{0.3} with $\R(A)=k$, then $M$ admits a formal,
CR transversal embedding into $\5 H^{2N+1}_{\l'}$ with $N=n+k$ and
$\l'=\min(\l+\s(A), N-\l-\s(A))$. Note that, since $\mathbb
H^{2N+1}_{\l'}$ is equivalent to $\mathbb H^{2N+1}_{N-\l'}$, it is
enough to show that there is an embedding into $\mathbb
H^{2N+1}_{\l'}$ with $\l'=\l+\s(A)$.

As pointed out in the introduction, if
$\R(A)=k$ and $\s(A)=s$, then there are formal power series
$(\phi_j(z,w))_{1\le j\le k}$ such that
$$
A(z,\bar z,w,\bar w)=
-\sum_{j=1}^{s} |\phi_j(z,w)|^2+\sum_{j=s+1}^{k} |\phi_j(z,w)|^2.
$$
The map $H(z,w):=(z,\phi(z,w),w)$, where
$\phi:=(\phi_1,\ldots,\phi_{k})$, defines a formal, CR transversal
embedding of $M$ into $\5 H^{2N+1}_{\l'}$ with $N=n+k$ and
$\l'=\l+s$. This completes the proof of Proposition~\ref{p0.2}.
\epf

\bpf[Proofs of Theorems $\ref{t0.4}$ and $\ref{t.extra}$] Let $M$,
$H_1$, $H_2$, $k_1$ and $k_2$ be as in Theorem~\ref{t0.4} and set
$N_1:=n+k_1$, $N_2:=n+k_2$. Observe that, as was mentioned in the
introduction, the embeddings $H^1$, $H^2$ must be CR transversal,
since $\l'=\l<n$. Choose a local coordinate system $(z,w)\in
\bC^n\times \bC$, vanishing at $p\in M$, such that $M$ is given by
an equation of the form \eqref{0.3}. Let $\sigma_q=\pm 1$ so that,
in the notation of Lemma \ref{l4.1}, \eqref{sigma} holds for
$H_q$. In view of the remark following the lemma, we must have
$\sigma_q=+1$ unless $\l=\l'=n/2$. By Lemma~\ref{l4.1}, there are
$T_q\in \Aut(\5 H^{2N_q+1}_\l,0)$, $q=1,2$, with the following
property: If we set $\tilde H_q:=T_q^{-1}\circ H_q$, then
$$\tilde H_q=(F_q,\phi_q,G_q)\in \C^n\times\C^{N_q-n}\times \C,$$
where $\phi_q$ have no constant or linear terms, and
$(F_q(z,w),G_q(z,w)=(z+f_q(z,w),\sigma_q w+g_q(z,w))$, where
$(f_q,g_q)$ both satisfy the normalization conditions \eqref{1.1}.

Let us first consider the case $\l=\l'<n/2$, so that $\sigma_q=+1$.
By making the local changes of coordinates
$(z_q,w_q)=(F_q(z,w),G_q(z,w))=(z+f_q(z,w),\sigma_q w+g_q(z,w))$,
for $q=1,2$, we observe, in view of the discussion in the proof of
Proposition~\ref{p0.2} above, that $M$ is defined, in the
coordinates $(z_q,w_q)$, by an equation of the form \eqref{0.3}
with $A=A_q$ given by
\begin{equation*}%\Label{4.4}
A_q(z,\bar z,w,\bar w):= \sum_{j=1}^{k_q} |\tilde
\phi_{q,j}(z,w)|^2, \quad \tilde \phi_q:=\phi_q\circ
(F_q,G_q)^{-1},
\end{equation*}
where, for convenience, we have dropped the subscript $q$ on the
variables. Since $k_1+k_2< n$ by the assumptions,
Theorem~\ref{t3.2} (or Theorem 0.3 for that matter) implies that
$(F_1,G_1)\equiv (F_2,G_2)$ (since $(F_1,G_1)\circ
(F_2,G_2)^{-1}:=\id+(\tilde f,\tilde g)$ where $(\tilde f,\tilde
g)$ also satisfies the normalization conditions \eqref{1.1}) and
\begin{equation}\Label{4.5}
\sum_{j=1}^{k_1} |\tilde\phi_{1,j}(z,w)|^2
\equiv \sum_{j=1}^{k_2} |\tilde\phi_{2,j}(z,w)|^2.
\end{equation}
%The last identity means that $(\tilde\phi_{1,j})_{1\le j\le k_1}$
%and $(\tilde\phi_{2,j})_{1\le j\le k_2}$ span the same space
%$V_{A_1}=V_{A_2}$ as defined in section \ref{s0}.
Let us write
$\hat\phi_1:=(\tilde\phi_1,0)\in\C^{k_1}\times\C^{k_2-k_1}$, so
that $\hat\phi_1$ has $k_2$ components (recall that $k_2\geq
k_1$). Then it follows from \eqref{4.5} and \cite[Proposition 3, pp 102]{Da}
that there exists a constant unitary transformation $U$
of $\bC^{k_2}$ such that $\hat\phi_1\equiv U\tilde \phi_2$.
Hence, by composing $\tilde H_2$ with a (unitary linear)
automorphism $T\in \Aut(\5 H^{2N_2+1}_\l,0)$, we obtain
$T\circ\tilde H_2 \equiv L\circ\tilde H_1$, where $L$ denotes the
linear embedding as the statement of Theorem~\ref{t0.4} (so that
$L\circ \tilde H_1=(F_1,\hat\phi_1,G_1)$.) The conclusion of
Theorem~\ref{t0.4}, with $\l'=\l<n/2$, follows from the
construction of $\tilde H_q$, $q=1,2$, and the easily verified
fact that any embedding $L\circ T_1$, with $T_1\in \Aut(\5
H^{2N_1+1}_\l,0)$, is of the form $T_2\circ L$ for some $T_2\in
\Aut(\5 H^{2N_2+1}_\l,0)$.

In the case $\l'=\l=n/2$, we cannot exclude the case
$\sigma_q=-1$. We consider the local coordinates $(z_q,
w_q)=(F_q(z,w),\sigma_q G_q(z,w))$ and proceed as above. The exact
same arguments show that $(F_1,\sigma_1 G_1)\equiv (F_2,\sigma_2
G_2)$ and that
\begin{equation}\Label{4.6}
\sigma_1\sum_{j=1}^{k_1} |\tilde\phi_{1,j}(z,w)|^2 \equiv\sigma_2
\sum_{j=1}^{k_2} |\tilde\phi_{2,j}(z,w)|^2,
\end{equation}
where $\tilde \phi_q=\phi\circ (F_q,\sigma_q G_q)^{-1}$). Now,
there are two cases to consider, namely $\sigma_1\sigma_2=1$ and
$\sigma_1\sigma_2=-1$. In the former case, the conclusion that
$H_2=T\circ L\circ H_1$, where $T$ and $L$ are as in the statement
of the theorem, follows as above.

In the latter case $\sigma_1\sigma_2=-1$, we conclude from
\eqref{4.6} that $\phi_{q,j}\equiv 0$ for $q=1,2$ and $j=1,\ldots,
N-n$. Hence, $M$, in the coordinates say $(z_1,w_1)$, is equal to
the quadric $\mathbb H^{2n+1}_{n/2}$ and, hence, this case never
occurs in the situation of Theorem \ref{t0.4}. This completes the
proof of Theorem \ref{t0.4}.

To prove Theorem \ref{t.extra}, we let $H_2=H$ and define $\tilde
H_1$ as above. Then, we define $H_1\colon (\mathbb
H^{2n+1}_{n/2},0)\to (\mathbb H^{2n+1}_{n/2},0)$, i.e. with
$k_1=0$, according to the sign of $\sigma_2$ as follows
\begin{equation}
H_1(z,w)=\left\{\begin{aligned} &(z,w),\quad \text{{\rm if
$\sigma_2=1$}}\\
&(z,-w),\quad \text{{\rm if $\sigma_2=-1$.}}\end{aligned}\right.
\end{equation}
Now, $\sigma_1\sigma_2=1$ and by proceeding as above, we conclude
that $H_2=T\circ L\circ H_1$, where $T$ and $L$ are as in Theorem
\ref{t.extra}. Recall that when $\sigma_2=-1$ we have renumbered
the coordinates $z'\in \bC^{N_2}$ so that \eqref{4.1} holds. If we
undo this reordering, then $L\circ H_2$ is equal to $L_-$ as
defined by \eqref{L-}. This completes the proof of Theorem
\ref{t.extra}. \epf

Let us briefly consider the case where $M$ is embedded into $\5
H^{2N+1}_{\l'}$ with $\l'>\l$. As demonstrated by
Example~\ref{e0.5}, we cannot hope for the full conclusion of
Theorem~\ref{t0.4} in this case. In what follows, we assume that
the signature $\l$ of $M$ is fixed and that there is a CR
transversal embedding of $M$ into $\5 H^{2(n+k_1)+1}_{\l'}$. For
simplicity, we shall only consider the case where $\l'<n-\l$, so
that $\sigma$, defined by \eqref{sigma}, must equal $+1$. (When
$\l'\geq n-\l$ there are different cases to consider, as in the
proofs of Theorems \ref{t0.4} and \ref{t.extra} above.) We shall
assume that the Hermitian product $\<\cdot,\cdot\>_{\l'}$ on
$\bC^{N}$, with $N:=n+k$, is given by \eqref{4.1}. A similar
argument to the one in the proof of Theorem~\ref{t0.4} above,
yields the following result; the details of modifying the argument
are left to the reader.

\bt\Label{t4.2} Let $M$ be a formal Levi nondegenerate
hypersurface in  $\bC^{n+1}$  of signature $\l$ at a point $p$.
Suppose that there are two formal, CR transversal embeddings
$H_q\: (M,p)\to (\5 H^{2(n+k_q)+1}_{\l_q},0)$, $q=1,2$, $($with
$k_2\geq k_1\ge 0$ and $\l\leq \l_q< n-\l)$. If $k_1+k_2<n$, then
there exist automorphisms $T_q\in \Aut(\5
H^{2(n+k_q)+1}_{\l_q},0)$
 and formal holomorphic coordinates $(z,w)\in \bC^n\times\bC$,
vanishing at $p\in M$, such that
\begin{equation}\Label{id-fin}
(T_q\circ
H_q)(z,w)=(z,\phi_q(z,w),w)\in\C^n\times\C^{k_q}\times\C, \quad
\<\phi_1,\phi_1\>_{\l_1-\l} \equiv \<\phi_2,\phi_2\>_{\l_2-\l},
\end{equation}
where the $\phi_q$'s have no constant or linear terms. \et

\br In the setting of Theorem~\ref{t4.2}, let us suppose that
$k_2-k_1\geq 2(\l_2-\l_1)$ (the other case can be treated
analogously). We write $\phi_q=(\phi_q^1,\phi_q^2)\in
\C^{\l_q-\l}\times \C^{k_q-\l_q+\l}$. It then follows from
Theorem~\ref{t4.2}, as in the  proof of Theorem~\ref{t0.4} above,
that there exists a unitary transformation $U$ of
$\bC^{k_2+\l_1-\l_2}$ such that
 $L(\phi_2^1,\phi_1^2)=U(\phi_1^1,\phi_2^2)$,
%or $(\phi_1^1,\phi_2^1)=U(\phi_1^2,\phi_2^2)$
%depending on the sign in \eqref{id-fin},
where $L$ is the linear embedding of $\bC^{k_1+\l_2-\l_1}$ into
$\bC^{k_2+\l_1-\l_2}$ via $z\mapsto (z,0)$.
%in the special case $k_2-k_1=2(\l_1-\l_2)$ (so that $k_2+\l_1-\l_2=k_1+\l_2-\l_1$), the
%conclusion is that $(\psi^2,\phi^1)=U(\psi^1,\phi^2)$.
For instance, in the situation of Example~\ref{e0.5} with $n \geq
2 $, where $k_1=0$, $\l_1=\l=0$, $k_2=2$, $\l_2=\l+1=1$ (so that
$\l_2<n-\l=n$ and $k_2+\l_1-\l_2=k_1+\l_2-\l_1=1$), we conclude
that there is a unimodular complex number $u$ such that
$\phi^1_2\equiv u\phi^2_2$. \er

\end{document}